# Converting of algebraic Diophantine equations to a diagonal form with the help of an integer non-orthogonal transformation, maintaining the asymptotic behavior of the number of its integer solutions

VICTOR VOLFSON

ABSTRACT. The author showed that any homogeneous algebraic Diophantine equation of the second order can be converted to a diagonal form using an integer non-orthogonal transformation maintaining asymptotic behavior of the number of its integer solutions. In this paper, we consider the transformation to the diagonal form of a wider class of algebraic second-order Diophantine equations, and also we consider the conditions for converting higher order algebraic Diophantine equations to this form. The author found an asymptotic estimate for the number of integer solutions of the diagonal Thue equation of odd degree with an amount of variables greater than two, and also he got and asymptotic estimates of the number of integer solutions of other types of diagonal algebraic Diophantine equations.

1. INTRODUCTION

This article is devoted to applications of algebra - the theory of affine transformations.

We study the number of integer solutions of algebraic Diophantine equations with integer coefficients in the paper. Asymptotic estimates for the number of integer solutions of Diophantine equations in [1-4], [6], [7] were carried out mainly for diagonal equations that correspond to the canonical equations of surfaces. Therefore, it is interesting to determine the class of transformations of an algebraic Diophantine equation that can lead to diagonal form of the equation, maintaining the asymptotic behavior of the number of its integer solutions.

We investigated the asymptotic estimate of the number of integer solutions of the diagonal Thue equation of even degree with an amount of variables greater than two in [4]. An asymptotic estimate of the number of integer solutions of the diagonal Thue equation of odd degree will be determined for the same number of variables in this paper. We will also consider asymptotic estimates of the number of integer solutions of other types of diagonal algebraic Diophantine equations.

______________________________________________

We considered an integer generalized orthogonal transformation in [5], which converts the algebraic Diophantine equation to diagonal form, maintaining the asymptotic behavior of the number of its integer solutions. However, this transformation does not make it possible to convert any homogeneous algebraic Diophantine equation of the second order to this form. Therefore, it is interesting to investigate non-orthogonal transformation in this direction.

It will be shown in the paper that with the help of an integer non-orthogonal transformation one can convert any homogeneous algebraic Diophantine equation of the second order to a diagonal form, maintaining the asymptotic behavior of the number of its integer solutions. We also consider the transformation of a wider class of algebraic second-order Diophantine equations corresponding to central and non-central hypersurfaces. The conditions for converting to the diagonal form of higher order algebraic Diophantine equations will be also considered.

## 2. ASYMPTOTIC ESTIMATE OF THE NUMBER OF INTEGER SOLUTIONS OF SERTAIN TYPES OF DIAGONAL ALGEBRAIC DIOPHANTINE EQUATIONS

We consider the asymptotic estimate of the number of integer solutions of the diagonal algebraic Diophantine equation of the degree $n$ for $k$ variables:

$$a_1 x_1^n + ... + a_{k-1} x_{k-1}^n + 2a_k x_k = 0. \tag{2.1}$$

Equation (2.1) corresponds to a paraboloid with the value $n = 2$.

If the coefficients are: $a_1 > 0, ..., a_{k-1} > 0, a_k < 0$, then equation (2.1) corresponds to an elliptic paraboloid. If the coefficients: $a_i, (i = 1, ..., k-1)$ have different signs, then equation (2.1) corresponds to a hyperbolic paraboloid.

It was shown in [4] that Thue equation of the degree $n$ ($n$ is even) of variables $k$:

$$a_1 x_1^n + ... + a_k x_k^n + a_0 = 0 \tag{2.2}$$

has the following upper asymptotic estimate for the number of integer solutions in a hypercube with side $[-N, N]$:

$$R_k(N) \ll N^{k-m+\epsilon}, \qquad (2.3)$$

where $m$ is the number of terms of equation (2.2) for which values $a_i > 0, (i = 1,...,k)$.

We call attention to the fact that equation (2.1) for a fixed value $x_k$ is Thue equation of order $n$ ($n$ is order) for $k-1$ variables. Therefore, based on (2.3), the following asymptotic upper bound is valid for the number of its integer solutions:

$$R_{k-1}(N) \ll N^{k-1-m+\epsilon}. \qquad (2.4)$$

Having in mind the estimate (2.4) and the fact that $x_k$ takes $2N+1$ integer values on the interval $[-N, N]$, we obtain the following asymptotic upper bound for the number of integer solutions of equation (2.1) for an even value $n$:

$$R_k(N) < O(N^{k-1-m+\epsilon}) \cdot O(N) = O(N^{k-m+\epsilon}) \text{ or } R_k(N) \ll N^{k-m+\epsilon}. \qquad (2.5)$$

It is known [6] that the number of integer solutions of an algebraic Diophantine equation (not necessarily diagonal) of the order $n$ of variables $k$ in a hypercube with side $[-N, N]$ has the following upper bound:

$$R_k(N) \ll N^{k-1+\epsilon}, \qquad (2.6)$$

where $\epsilon$ is an arbitrarily small positive number.

It is also known Pila formula [7] for an irreducible algebraic Diophantine equation (not necessarily diagonal) of the order $n (n > 2)$ of variables $k$ ($k > 2$) for the number of integer solutions in a hypercube with side $[-N, N]$:

$$R_k(N) \ll N^{k-2+\epsilon+1/n}, \qquad (2.7)$$

compare with (2.6).

Now we give an upper bound for the number of integer solutions of the diagonal Diophantine Thue equation of the order $n (n > 2)$ (odd) of variables $k$ ($k > 2$):

$$a_1 x_1^n + ... + a_k x_k^n + a_0 = 0, \qquad (2.8)$$

where all the coefficients $a_i$ are integers, distinct from $0$.

If the value $k = 2$, then we'll obtain the Thue equation from two variables, which will have a finite number of solutions for $n > 2$, that is:

$$R_2(N) = O(1). \tag{2.9}$$

Based on (2.9), we'll get on an interval $[-N, N]$:

$$R_k(N) \leq O(1) \cdot O(N^{k-2}) = O(N^{k-2}) \text{ or } R_k(N) \ll N^{k-2+\epsilon}. \tag{2.10}$$

if all other variables in equation (2.8) take integer values.

We note that the estimate (2.10) is more precise than (2.7).

Now we consider the diagonal algebraic Diophantine equation of the order $n(n > 2)$ (odd) for $k$ variables:

$$a_1 x_1^n + \ldots + a_{k-1} x_{k-1}^n + 2a_k x_k = 0. \tag{2.11}$$

We'll obtain Thue equation of order $k-1$ for which (2.10) is satisfied:

$$R_{k-1} \ll N^{k-3+\epsilon}, \tag{2.12}$$

if the value $x_k$ is fixed.

Having in mind (2.12), we'll obtain an upper bound for the number of integer solutions of equation (2.11) on the interval $[-N, N]$:

$$R_k(N) \ll N^{k-2+\epsilon}, \tag{2.13}$$

if $x_k$ takes integer values.

One can obtain a better upper bound for the number of integer solutions for the Diophantine equations (2.8) and (2.11).

Diophantine equation (2.8) does not have solutions in integers in the case when all the coefficients $a_i > 0 (i = 0, \ldots, k)$ and $x_i > 0 (i = 1, \ldots, k)$. Equation (2.8) has a finite number of solutions in integers if $a_i > 0 (i = 1, \ldots, k)$ and $x_i > 0 (i = 1, \ldots, k)$. Thus, equation (2.8) has solutions in integers in both considered cases.

We'll replace equation (2.8) by the equation:

$$-(a_1 x_1^n + \ldots + a_k x_k^n + a_0) = 0, \tag{2.14}$$

in which more than half of the variables will have nonnegative values, if more than half values of variables of equation (2.8) $x_i < 0$.

Thus, the number of integer solutions of equation (2.8) does not exceed $O(N^{k/2})$. Consequently, the number of integer solutions of equation (2.8), depending on the ratio of the values of the variables $x_i$, takes values from $O(N^0)$ to $O(N^{k/2})$. Hence we obtain an upper bound for the number of integer solutions of equation (2.8):

$$R_k(N) << N^{k/2+\epsilon}. \qquad (2.15)$$

An analogous upper bound for (2.15) holds for the number of integer solutions of equation (2.11).

The formula (2.15) is valid for even values $k$. The following formula holds for odd values $k$:

$$R_k(N) << N^{[k/2]+1+\epsilon}. \qquad (2.16)$$

The formula (2.15) is more exact than formulas (2.10) and (2.13) for $k \geq 6$, and formula (2.16) is more exact than (2.10) and (2.13) for $k \geq 7$. It is necessary to use formulas (2.10) and (2.13) for smaller values $k$.

For example, we consider the following Diophantine Thue equation:

$$x_1^3 + x_2^3 + x_3^3 - 1 = 0. \qquad (2.17)$$

Therefore, having in mind (2.10), we obtain the following upper bound for the number of integer solutions of equation (2.17) in this case when the value $k = 3$:

$$R_3(N) << N^{1+\epsilon}. \qquad (2.18)$$

Based on [7] this coincides with the upper bound for the number of integer solutions of equation (2.17).

We show that the estimate (2.18) can't be improved.

We represent equation (2.17) in the form: $f = f_1 + f_2$, where $f_1 = x_1^3 + x_2^3 = 0$ and $f_2 = x_3^3 - 1 = 0$.

The homogeneous equation $f_1 = 0$ has only integer solutions on the line $x_1 + x_2 = 0$, that is, it has $O(N)$ integer solutions in the square with a side $[-N, N]$.

The equation $f_2 = 0$ has only one integer solution $x_3 = 1$. Therefore, equation (2.17) has the following lower bound for the number of integer solutions in a cube with side $[-N, N]: R_3(N) \gg N$. Compare this with the upper bound (2.18).

## 3. CONVERSION OF THE ALGEBRAIC DIOPHANTINE EQUATION TO THE DIAGONAL FORM WIT THE HELP OF AN INTEGER (NOT NECESSARILY) ORTOGONAL TRANSFORMATION MAINTING THE ASYMPTOTIC OF THE NUMBER OF ITS INTEGER SOLUTIONS

It was shown in [5] that the asymptotic of the number of integer solutions of an algebraic Diophantine equation is preserved for an integer affine transformation in a hypercube with side $[-N, N]$. Therefore, our goal will be to obtain an integer affine transformation that leads the algebraic Diophantine equation to the diagonal form.

Assertion 1

There is a linear transformation (possibly not orthogonal), with rational coefficients, converting quadratic form with rational coefficients to canonical form with rational coefficients.

Proof

It is known that one can bring any quadratic form to the canonical form by means of a Lagrange transformation, which is linear. On the other hand, all transformations with the help of Lagrange method take place in the field, so the transformation will have rational coefficients if coefficients of the form are rational numbers. Therefore, the canonical form also has rational coefficients.

Corollary 1

There exists an integer linear (not necessarily orthogonal) transformation that gives a homogeneous or second order Thue equation to the diagonal form, maintaining the asymptotic of their integer solutions in a hypercube with a side $[-N, N]$.

Proof

There exists a linear (not necessarily orthogonal) transformation with rational coefficients, which gives a quadratic form with integer coefficients to the canonical form with rational coefficients (based on Assertion 1). After it we'll perform a homothetic transformation with a coefficient equal to the smallest common denominator of the rational coefficients of the first transformation. We obtain as a result an integer linear (not necessarily orthogonal) transformation that leads to a diagonal form of a homogeneous equation or a Thue equation of the second order.

Any integral affine transformation preserves the asymptotic of integer solutions of an algebraic Diophantine equation [5] in a hypercube with side $[-N, N]$, therefore, the above resulting integer linear transformation also preserves this asymptotic behavior.

Corollary 2

If the conversion involving the Lagrange transformation and the transfer:

$$x_1 = c_{11} x'_1 + ... + c_{1k} x'_k + c_1, ... x_k = c_{k1} x'_1 + ... + c_{kk} x'_k + c_k,$$

which has integer values of the coefficients $c_1, ... c_k$, leads an algebraic Diophantine equation of the second order to the diagonal form, then with the help of an integer affine (not necessarily orthogonal) transformation, then this Diophantine equation can be reduced to a diagonal form maintaining the asymptotic of integer solutions of the original equation in the hypercube with side $[-N, N]$.

Proof

Suppose that a Diophantine equation of the second order:

$$x_1 = c_{11} x'_1 + ... + c_{1k} x'_k + c_1, ... x_k = c_{k1} x'_1 + ... + c_{kk} x'_k + c_k, \quad (3.1)$$

where $c_1, ... c_k$ are integers, is converted to a diagonal form by means of the transformation involving the Lagrange transformation and parallel transfer.

Based on Assertion 1, all the coefficients of the transformation (3.1) $c_{ij}$ are rational numbers, so after the transformation of the deformation:

$$x'_1 = t_1 x''_1, ..., x'_k = t_k x''_k, \quad (3.2)$$

where $t_j$ are the least common multiple of the denominator values $c_{ij}$ ($j = 1,...,k$) is taken in (3.2), we obtain the resulting integer affine transformation that converts the original equation to diagonal form with integer coefficients.

Having in mind [5], the diagonal algebraic Diophantine equation obtained after this integer transformation has the same asymptotic of integer solutions in a hypercube with side $[-N, N]$ as the original non diagonal algebraic Diophantine equation.

Let us explain the above with an example.

It is necessary to reduce the Diophantine equation to the diagonal form:

$$4x_1^2 + 9x_2^2 + 12x_1x_2 + 8x_1 + 2x_2 + 24 = 0. \quad (3.3)$$

We write equation (3.3) in the form:

$$4(x_1 + 3x_2/2)^2 + 8x_1 + 2x_2 + 24 = 0.$$

Let us make the Lagrange transformation:

$$x_1 = x_1' - 3x_2'/2, \; x_2 = x_2'. \quad (3.4)$$

We transform (3.3) with the help (3.4) and obtain:

$$4x_1'^2 + 8x_1' - 10x_2' + 24 = 0. \quad (3.5)$$

We write equation (3.5) in the form:

$$4(x_1' + 1)^2 - 10(x_2' - 2) = 0. \quad (3.6)$$

Let's make the transformation of the origin transfer for equation (3.6):

$$x_1' = x_1'' - 1, \; x_2' = x_2'' + 2. \quad (3.7)$$

We obtain the equation after performing the transformation (3.7):

$$4x_1''^2 - 10x_2'' = 0. \quad (3.8)$$

We obtain the total transformation after the transformation (3.7) which is not an integer, but it has integer free coefficients:

$$x_1 = x_1'' - 3x_2''/2 - 4, \; x_2 = x_2'' + 2. \quad (3.9)$$

Let's do the deformation transformation to get the integer transformation from (3.9):

$$x_1'' = x_1''', x_2'' = 2x_2'''. \tag{3.10}$$

We obtain the resulting integer affine transformation after transformation of the deformation (3.10):

$$x_1 = x_1''' - 3x_2''' - 4, \quad x_2 = 2x_2''' + 2. \tag{3.11}$$

The integer resulting affine transformation (3.11) converted the original equation (3.3) to the diagonal form:

$$x_1'''^2 - 5x_2''' = 0. \tag{3.12}$$

The diagonal algebraic Diophantine equation (3.12) obtained after the resulting integer transformation (3.11) has the same asymptotic of integer solutions as the original non-diagonal algebraic Diophantine equation (3.3) in the hypercube with side $[-N, N]$.

It is well known that the affine transformation preserves the parallelism of hyper-lines and hyperplanes, therefore the hypercube in the general case is transformed into an inclined hyperparallelulepiped under an affine transformation. The hypercube goes to a direct hyperparallelepiped under a generalized orthogonal transformation preserving angles and the hypercube is transformed into the same hypercube under an orthogonal transformation preserving angles and distances.

Based on properties of an affine transformation, the volume of an oblique hyperparallelepiped (obtained after an affine transformation with a matrix $C$) of $n$-dimensional asymptotic hypercube with side $[-N, N]$ is equal to:

$$V_p = (2N)^n \cdot |det(C)|. \tag{3.13}$$

All integer points in a hypercube with a side $[-N, N]$ (under an integral affine transformation) go to the corresponding integer points in the inclined hyperparallelepiped obtained after the transformation.

Therefore, starting from (3.13) if $|det(C)| > 1$, then the density of integer solutions of the Diophantine equation decreases (the density of integer solutions of the Diophantine equation is the quotient of the number of integer solutions of the Diophantine equation and the volume of the hypercube).

It is known that a unimodular transformation is a linear integer transformation with a matrix $U$ for which $|det(U)|=1$.

Two statements follow from the above.

Assertion 2

The density of integer solutions of the Diophantine equation is preserved for unimodular transformation.

The number of integer solutions of the Diophantine equation in a hypercube is equal to the number of integer solutions of the Diophantine equation obtained after the transformation, which are in obtained oblique hyperparallelepiped.

Corollary

The transformation will be unimodular if it is an integer Lagrange transformation and all its corner minors are non-zero. Therefore, assertion 2 applies to it.

On the other hand, it converts the quadratic form of the Diophantine equation to the canonical form:

$$M_1 x_1^2 + (M_2/M_1) x_2^2 + ... + (M_r/M_{r-1}) x_r^2, \qquad (3.14)$$

where $M_i$ is the angular minor of order $i$, and $r$ is the rank of the quadratic form.

Proof

The matrix of the Lagrange transform has an upper triangular form with ones on the diagonal, since all the angular minors of the transformation matrix are nonzero.

Therefore, the determinant of such a matrix it is 1 and is unimodular, since the Lagrange transformation is linear and integer by hypothesis.

On the other hand, the condition of the Jacobi theorem is satisfied, since all the angular minors of the transformation matrix are nonzero. Jacobi's theorem asserts that the canonical form of a quadratic form is defined by formula (3.14) in this case.

Assertion 3

Any integer linear orthogonal transformation is unimodular and, therefore, the number of integer solutions of the Diophantine equation is preserved in the hypercube under a given transformation.

Proof

The module of the determinant of the transformation will be equal to one if this transformation is orthogonal.

This transformation is unimodular, since by assumption it is linear and integer.

Based on Assertion 2, the unimodular transformation preserves the number of integer solutions of the Diophantine equation in the hypercube.

The resulting hypercube contains the same number of integer solutions of the Diophantine equation as the original one, since the hypercube converts into the same hypercube under the orthogonal transformation.

The sines and cosines of the rotation angles around the coordinate axes take only the values 0 and 1, 0 and -1, respectively, for an integer linear orthogonal transformation. Therefore, the rotation angles around the coordinate axes are multiples $\pi/2$ for a given transformation, there is also possible axial symmetry relative to the coordinate axes or the composition of these transformations.

I will explain assertion 2 in the example.

There is a homogeneous Diophantine equation of the second order:

$$x_1^2 + 2x_1 x_2 - 3x_2^2 = 0. \tag{3.15}$$

Let us single out the complete square in the quadratic form of equation (3.15):

$$(x_1 + x_2)^2 - x_2^2 - 3x_2^2 = (x_1 + x_2)^2 - 4x_2^2 = (x_1 + 3x_2)(x_1 - x_2) = 0. \tag{3.16}$$

Integer solutions of equation (3.16) are on the line $x_2 = -x_1/3$, so the equation has such number of integer solutions: $R_2(N) = 1 + 2N/3$ in the square with the side $[-N, N]$. The equation also has integer solutions on the line $x_2 = x_1$ in the square with the side $[-N, N]$ - $R_2(N) = 2N + 1$.

The total number of integer solutions of equation (3.15) is equal to:

$$R_3(N) = 2N + 2N/3 + 1 = 1 + 8N/3, \qquad (3.17)$$

since both equations have a trivial solution.

Having in mind (3.17), the number of integer solutions of equation (3.15) in a square with side $[-3, 3]$ is equal to:

$$R_3(3) = 9. \qquad (3.18)$$

We make the change of variables in equation (3.15) by the Lagrange method:

$x_1' = x_1 + x_2, x_2' = x_2$ , which corresponds to the transformation:

$$x_1 = x_1' - x_2', x_2 = x_2'. \qquad (3.19)$$

Based on of (3.19), the matrix of the transformation by the Lagrange method has the form:

$$C = \begin{pmatrix} 1 & -1 \\ 0 & 1 \end{pmatrix}, \qquad (3.20)$$

those $det(C) = 1$.

Therefore, the integer matrix (3.20) is unimodular and the conditions of Assertion 2 are satisfied.

We obtain a homogeneous Diophantine equation with the canonical quadratic form from (3.15) after the unimodular transformation (3.19):

$$(x_1')^2 - 4(x_2')^2 = (x_1' + 2x_2')(x_1' - x_2') = 0. \qquad (3.21)$$

Equation (3.21) has integer solutions. These are on straight lines:

$$x_2' = -x_1'/2, x_2' = x_1'. \qquad (3.22)$$

The square with a side $[-3, 3]$ under transformation converts into a parallelogram with vertices:

$$(0, 3), (6, 3), (0, -3), (-6, -3). \qquad (3.23)$$

Equation (3.21) in the parallelogram (3.23) on the lines (3.22) also has 9 integer solutions, which corresponds to Assertion 2.

Suppose we have a homogeneous nondiagonal algebraic Diophantine equation of the second order for three variables:

$$a_{11}x_1^2 + a_{22}x_2^2 + a_{33}x_3^2 + 2a_{12}x_1x_2 + 2a_{13}x_1x_3 + 2a_{23}x_2x_3 = 0. \tag{3.24}$$

Let $M_i$ is an angular minor of order $i$ (not equal to 0) ($i = 1, 2, 3$) of a matrix of quadratic form:

$$A_3 = \begin{pmatrix} a_{11} & a_{12} & a_{13} \\ a_{21} & a_{22} & a_{23} \\ a_{31} & a_{32} & a_{33} \end{pmatrix}. \tag{3.25}$$

It is necessary to find the corresponding homogeneous diagonal algebraic Diophantine equation based on coefficients of the nondiagonal equation.

It is easy to show that the matrix of the Lagrange transformation for equation (3.24) will have the form:

$$C_3 = \begin{pmatrix} 1 & -a_{12}/M_1 & -a_{12}(a_{11}a_{23} - a_{12}a_{13})/M_1M_2 \\ 0 & 1 & (a_{12}a_{13} - a_{11}a_{23})/M_2 \\ 0 & 0 & 1 \end{pmatrix}. \tag{3.26}$$

if $M_i$ is not equal to 0.

Based on Jacobi's theorem (after the transformation by the Lagrange method), equation (3.24) takes the form:

$$M_1(x_1')^2 + (M_2/M_1)(x_2')^2 + (M_3/M_2)(x_3')^2 = 0. \tag{3.27}$$

The matrix (3.26) will be an integer (unimodular, since $det(C_3) = 1$) if $a_{12}$ and $a_{12}a_{13} - a_{11}a_{23}$ are integers, multiples $M_1$, and at the same time $a_{12}(a_{11}a_{23} - a_{12}a_{13})$ are integers, multiples $M_1M_2$.

The transformation matrix by the Lagrange method (3.26) is not integer in the remaining cases, and it is necessary an additional deformation transformation to bring the equation (3.24) to diagonal form:

$$x_1' = k_1 x_1'', x_2' = k_2 x_2'', x_3' = k_3 x_3''. \tag{3.28}$$

Equation (3.27) takes the diagonal form, taking into account the deformation transformation (3.28):

$$M_1 k_1^2 (x_1'')^2 + M_2 k_2^2 / M_1 (x_2'')^2 + M_3 k_3^2 / M_2 (x_{3''})^2 = 0. \tag{3.29}$$

Values of the stretching coefficients $k_1, k_2, k_3$ in equation (3.29) are equal to:

1. If $a_{12}$ and $a_{12}(a_{11}a_{23} - a_{12}a_{13})$ are integers, then: $k_1 = 1, k_2 = M_1, k_3 = M_1 M_2$.

2. If $a_{12}$ is an integer, but $a_{11}a_{23} - a_{12}a_{13}$ is fractional one, then: $k_1 = 1, k_2 = M_1, k_3 = 2M_1 M_2$.

3. If $a_{12}$ is a fractional, but $a_{11}a_{23} - a_{12}a_{13}$ is integer one, then: $k_1 = 1, k_2 = 2M_1, k_3 = M_1 M_2$.

4. If $a_{12}$ is a fractional, but $a_{11}a_{23} - a_{12}a_{13}$ is integer one, and $a_{12}(a_{11}a_{23} - a_{12}a_{13})$ is fractional one, then: $k_1 = 1, k_2 = 2M_1, k_3 = M_1 M_2$.

5. If $a_{12}$ is a fractional and $a_{11}a_{23} - a_{12}a_{13}$ is fractional one, then: $k_1 = 1, k_2 = 2M_1, k_3 = M_1 M_2$. (3.30)

Thus, having in mind (3.30), we can conclude that already for 3 variables the number of consideration cases is large, especially for a larger number of variables. Therefore, finding the coefficients of the diagonal equation without finding the transformation matrix for a large number of variables is laborious.

We consider the particular case of the unimodular matrix (3.26), where $|M_1| = |M_2| = |M_3| = 1$, where $M_i$ is the angular minor order $i$ ($i = 1, 2, 3$) of the matrix of the quadratic form (3.25).

Equation (3.27) takes the normal form in this case:

1. If $M_1 = -1, M_2 = -1, M_3 = -1$, then we'll get Fermat equation $-x_1^2 + x_2^2 + x_3^2 = 0$.

2. If $M_1 = -1, M_2 = -1, M_3 = 1$, then we'll get Fermat equation $-x_1^2 + x_2^2 - x_3^2 = 0$.

3. If $M_1 = -1, M_2 = 1, M_3 = -1$, then we'll get the equation only with a trivial solution $-x_1^2 - x_2^2 - x_3^2 = 0$.

4. If $M_1 = -1, M_2 = 1, M_3 = 1$, then we'll get Fermat equation $-x_1^2 - x_2^2 + x_3^2 = 0$.

5. If $M_1 = 1, M_2 = -1, M_3 = -1$, then we'll get Fermat equation $x_1^2 - x_2^2 + x_3^2 = 0$.

6. If $M_1 = 1, M_2 = -1, M_3 = 1$, then we'll get Fermat equation $x_1^2 - x_2^2 - x_3^2 = 0$.

7. If $M_1 = 1, M_2 = 1, M_3 = -1$, then we'll get Fermat equation $x_1^2 + x_2^2 - x_3^2 = 0$.

8. If $M_1 = 1, M_2 = 1, M_3 = 1$, then we'll get the equation only with a trivial solution $x_1^2 + x_2^2 + x_3^2 = 0$. (3.31)

An asymptotic upper and lower bound for the number of integer solutions of the second-order Fermat equation was obtained in (5) for the cube with side $[-N, N]$:

$$N \ll R_3(N) \ll N \ln(N). \tag{3.32}$$

Consequently, the asymptotic estimate (3.32) is also valid for all nondiagonal algebraic Diophantine equations of second-order for 3 variables with $|M_1| = |M_2| = |M_3| = 1$ in the cube with side $[-N, N]$, except for cases 3 and 8 in (3.31).

Based on of Jacobi's theorem, we obtain a diagonal equation in the normal form after the Lagrange transformation for a Diophantine equation of the second order for $n$ variables when the condition $|M_1| = ... = |M_n| = 1$ and integrity of the transformation matrix performed:

$$x_1^2 + ... x_r^2 - x_{r+1}^2 - ... - x_n^2 = 0, \tag{3.33}$$

where $r$ is a positive and $n = r$ is negative form index.

If $n = r$, then equation (3.33) will have only one integer trivial solution. If $n - r > 0$, then equation (3.33) will have an infinite number of integer solutions.

Assertion 4

Equation (3.33) will have the following lower bound for the number of integer solutions in a hypercube with side $[-N, N]$, if $n - r > 0$:

$$R_n(N) \geq (4N+1)^{n-r}. \tag{3.34}$$

We prove this fact by mathematical induction on the number of variables. The first step induction:

We consider the equation:

$$x_1^2 + ... x_r^2 - x_{r+1}^2 = 0. \tag{3.35}$$

Let's represent equation (3.35) as the sum of two equations: $f_1 = x_2^2 + ... x_{r-1}^2 = 0, f_2 = x_r^2 - x_{r+1}^2$.

The equation $f_1 = 0$ has one trivial solution. The equation $f_2 = (x_r + x_{r+1})(x_r - x_{r+1}) = 0$ has in the square with the side $[-N, N]$ - $4N+1$ integer solutions. Therefore, the equation $f = f_1 + f_2 = 0$ has the following lower bound for the number of integer solutions in the square with the side $[-N, N]$:

$$R_{r+1}(N) \geq 4N+1,$$

which corresponds to the Assertion 4.

We make the assumption (by induction) that the assertion (3.34) is satisfied for $n - r = i$ when $r \geq i+1$ (if the condition is not satisfied, then we put the minus sign before the equation), ie:

$$R_{r+i}(N) \geq (4N+1)^i. \tag{3.36}$$

Let us show that the assertion is fulfilled for $n - r = i+1$ ($r \geq i+2$). If the condition does not hold, then we'll put a minus sign before the equation.

Equation (3.35) can be represented in this case in the form:

$f_3 = f_4 + f_5 = 0$, где $f_4 = x_1^2 + ... + x_{r-i-2}^2 + ... + x_{r-1}^2 - x_{r+i}^2 = 0$, $f_5 = x_{r-i-1}^2 - x_{r+i+1}^2 = 0$.

Based on the induction the estimate (3.36) holds for the equation - $f_4 = 0$. The following estimate $R_2(N) \geq 4N+1$ holds for the equation - $f_5 = 0$. Therefore, the following estimate $R_{r+i+1}(N) \geq (4N+1)^i(4N+1) = (4N+1)^{i+1}$ holds for the equation $f_3 = f_4 + f_5 = 0$, which corresponds to Assertion 4.

Based on the assertion, the estimate (3.34) is true for any nondiagonal Diophantine equation of second order for $n$ variables, when $|M_1| = ... = |M_n| = 1$ and the Lagrange transformation matrix $C_n$ is an integer.

Now we consider the conversion to diagonal form of algebraic Diophantine equations of higher degrees ($n > 2$). Forms of high orders are also converted to diagonal form in some cases [8].

Assertion 5

If there is a Lagrange transformation:

$$x_1' = c_{11}x_1 + ... + c_{1k}x_k, ..., x_k' = c_{k1}x_1 + ... + c_{kk}x_k, \quad (3.37)$$

which will convert a form to canonical species:

$$a_1'(x_1')^n + ... + a_k'(x_k')^n. \quad (3.38)$$

Then there is a deformation conversion:

$$x_j' = t_j x_j'' (j = 1,...,k) \quad (3.39)$$

where $t_j$ is the smallest common multiple $c_{ij}$, after which, having in mind (3.38), (3.39), the form of second order will be:

$$a_1' t_1^n (x_1'')^n + ... + a_k' t_k^n (x_k'')^n \quad (3.40)$$

with integer coefficients.

The proof of (3.40) follows from the fact that the Lagrange transformation (3.37) has coefficients $c_{ij}$ which (as was said above) are rational numbers in the case when the original form has integer coefficients.

Corollary 1

If the non-diagonal algebraic Diophantine equation of order $n$ (corresponding to the central surface for which the center is not at the origin) by means of an affine transformation containing the Lagrange transformation and the subsequent transfer:

$$x'_1 = c_{11}x_1 + ... + c_{1k}x_k + c_1, ..., x'_k = c_{k1}x_1 + ... + c_{kk}x_k + c_k,$$

where $c_1, ... c_k$ are integers, can be brought to the form:

$$a'_1(x'_1)^n + ... + a'_k(x'_k)^n + a'_0 = 0, \tag{3.41}$$

where $a'_0$ is an integer.

Then, using the conversion of the deformation:

$$x'_j = t_j x''_j (j = 1, ..., k), \tag{3.42}$$

where $t_{ij}$ is the smallest common multiple $c_{ij}$ (based on assertion 5) the equation:

$$a'_1 t_1^n (x''_1)^n + ... + a'_k t_k^n (x''_k)^n + a'_0 = 0 \tag{3.43}$$

will have integer coefficients.

Note

If the non-diagonal Diophantine equation corresponding to the central surface of order n which center is at the origin, can be reduced to the form (3.41) with the help of the Lagrange transformation (3.37). Then, using the transformation of the deformation (3.42) (based on Assertion 5) the equation (3.41) will can be converted to the form (3.43) with integer coefficients.

Equation (3.43) is a diagonal Thue equation of order n. An asymptotic estimate of the number of its integer solutions in a hypercube with a side was $[-N, N]$ considered earlier in this paper.

As was shown in [5], if there is an integer affine transformation that leads to a diagonal algebraic Diophantine equation, then the asymptotic estimates of the number of their integer solutions will coincide.

Therefore, the asymptotic estimate of the number of integer solutions of a non-diagonal Diophantine algebraic equation of order n corresponding to the central surface for the conditions indicated in Corollary 1 coincides with the asymptotic estimate of the number of integer solutions of equation (3.43) in a hypercube with a side $[-N, N]$.

Corollary 2

If an non-diagonal algebraic Diophantine equation of order n, corresponding to a non-central surface, by means of an affine transformation containing the Lagrange transformation and the subsequent transfer:

$$x'_1 = c_{11}x_1 + ... + c_{1k}x_k + c_1, ..., x'_k = c_{k1}x_1 + ... + c_{kk}x_k + c_k,$$

where are integers, can be brought to the form:

$$a'_1(x'_1)^n + ... + a'_{k-1}(x'_{k-1})^n + 2a'_k x'_k = 0.$$ (3.44)

Then, using the conversion of deformation:

$$x'_j = t_j x''_j (j = 1,...,k),$$

where $t_{ij}$ is the smallest common multiple $c_{ij}$, equation (3.44) will can be reduced to the diagonal form:

$$a'_1 t_1^n (x''_1)^n + ... + a'_{k-1} t_{k-1}^n (x''_{k-1})^n + 2a'_k t_k x'_k = 0$$ (3.45)

with integer coefficients.

The proof of Corollary 2 follows from the fact that all the coefficients of the Lagrange transform $c_{ij}$ are rational numbers.

An asymptotic estimate of the number of integer solutions of the diagonal algebraic Diophantine equation (3.45) in a hypercube with a side $[-N, N]$ was considered earlier in this paper.

As was shown above, if there is an integer affine transformation that converts algebraic Diophantine equation to a diagonal form, then the asymptotic estimates of the number of their integer solutions will be the same. Therefore, the asymptotic estimate of the number of integer solutions of a non-diagonal Diophantine algebraic equation of order n that corresponds to a non-central surface under the conditions indicated in Corollary 2 coincides with the asymptotic estimate of the number of integer solutions of equation (3.45) in a hypercube with a side $[-N, N]$.

Let's consider an illustrative example.

There is a nondiagonal algebraic Diophantine equation of the third order:

$$8x^3 + 3xy^2 + 12x^2 y + 28y^3 = 0. \qquad (3.46)$$

It is required to find a diagonal algebraic equation of the third order that has the same asymptotics of the number of integer solutions in a square with a side $[-N, N]$ as the nondiagonal equation (3.46), and It is required to determine the asymptotics of its integer solutions.

Decision

We single out the complete cube from the cubic form on the left of (3.46) by the Lagrange method:

$$8x^3 + 3xy^2 + 12x^2 y + 28y^3 = 8(x + y/2)^3 - y^3 + 28y^3 = 8(x + y/2)^3 + (3y)^3 = 0. \qquad (3.47)$$

We replace the variables in (3.47):

$$x' = x + y/2,\ y' = 3y$$

or we make a linear transformation:

$$x = x' - y'/6,\ y = y'/3. \qquad (3.48)$$

We obtain the equation:

$$8x'^3 + y'^3 = 0. \qquad (3.49)$$

performing the transformation (3.48) c (3.47).

The transformation (3.48) is not an integer, so we perform the transformation of the deformation:

$$x' = t_1 x'',\ y' = t_2 y'', \qquad (3.50)$$

where $t_1 = 1, t_2 = 6$ $t_1 = 1, t_2 = 6$.

Having in mind (3.48) we obtain the integer transformation:

$$x = x'' - y'', \quad y = 2y''. \tag{3.51}$$

after transformation of the deformation (3.50).

Based on the resultant integer transformation (3.51) we obtain an integer diagonal algebraic Diophantine equation from equation (3.46):

$$(2x'')^3 + (6y'')^3 = 0, \tag{3.52}$$

which has integer solutions on the line:

$$x'' + 3y'' = 0.$$

Thus, the diagonal Diophantine algebraic equation (3.52) has the asymptotics of integer solutions in the square with the side $[-N, N]$, which coincides with the asymptotics of the integer solutions of the nondiagonal Diophantine equation (3.46), which is equal to:

$$R_2(N) = O(N).$$

## 4. CONCLUSION AND SUGGESTIONS FOR FURTHER WORK

The next article will continue to study the asymptotic behavior of the number of integer solutions for the algebraic Diophantine equation.

## 5. ACKNOWLEDGEMENTS

Thanks to everyone who has contributed to the discussion of this paper.